\documentclass[11pt]{article}
\usepackage{amsmath, amsthm,amsfonts,amssymb,amscd}
\usepackage{color}
\usepackage{geometry}
\usepackage{tikz}
\usepackage{subfigure}
\numberwithin{equation}{section}
\usepackage[T1]{fontenc}
\usepackage{cite}
\usepackage{enumerate}

\usepackage{float}
\usepackage{soul}
\usepackage{marginnote}
\def\tto{\;{\lower 1pt \hbox{$\rightarrow$}}\kern -10pt
	\hbox{\raise 2pt \hbox{$\rightarrow$}}\;}

\def\h{\hfill\Box}

\def\h{\hfill\triangle}

\newcounter{lk}

\begin{document}
	
\newtheorem{Theorem}{Theorem}[section]
\newtheorem{Proposition}[Theorem]{Proposition}
\newtheorem{Remark}[Theorem]{Remark}
\newtheorem{Lemma}[Theorem]{Lemma}
\newtheorem{Corollary}[Theorem]{Corollary}
\newtheorem{Definition}[Theorem]{Definition}
\newtheorem{Example}[Theorem]{Example}
\newtheorem{Fact}[Theorem]{Fact}
\renewcommand{\theequation}{\thesection.\arabic{equation}}
	\normalsize
	\def\proof{
		\normalfont
		\medskip
		{\noindent\itshape Proof.\hspace*{6pt}\ignorespaces}}
	\def\endproof{$\h$ \vspace*{0.1in}}

	\title{\small \bf The extended eigenvalues of composition operator on Bergman space and Fock space}
	
	\date{}
	\author{  Jiale Li \footnote{Math, Chongqing University, Chongqing,  China; email: 20165726@cqu.edu.cn}}
	\maketitle
{\small \begin{abstract}
A complex scalar~$\lambda$~ is said to be an extended eigenvalue of a bounded linear operator A on a complex Hilbert space if there is a nonzero operator ~$X$~ such that ~$AX=\lambda XA$~. There are some solutions to the problem of computing the extended eigenvalues for composition operators ~$C_{\varphi}$~induced on the Fock space and Bergman space by linear fractional transformations~$\varphi$~of the complex plane in this paper.\\
\end{abstract}}
\section{introduction}
\subsection{total discussion}
\par Our research base on operator theory in function spaces, it is an important role under functional analysis which focus on discussing the properties of some specific linear operators on specific function space. Generally speaking, the function space we interact is an analytic function space on the complex plane, including Hardy space, Bergman space, Fock space, etc. Here are several books on these function space [1-3]. We usually define three types of operators on these function spaces: Toeplitz operators, Hankel operators and composition operators,those operators go through the basic theory. The common feature of these operators is that most of their properties are affected by their sign function. They not only have the common properties of linear operators in function spaces, but also own their special properties that are worth going deeply. Research in this field tend to study the boundedness, compactness, Fredholm property, reversibility, spectrum, etc. Particularly, whether these properties strongly correspond to some characteristics of their sign function, these characteristics of sign function and the corresponding relationship between sign function and operators are exactly what we addicted to discuss.
\par  Research on these operator theories is diverse and extremely rich, some of them even form unique deep level theoretical system. Reference [4] introduces the basic concept and theory of Toeplitz operators, Douglas provides an in-depth introduction to Toeplitz operator theory on Hardy spaces in Reference [5], and Reference [6] is a specialized work on Toeplitz operator theory. Reference [7] systematically introduces Hankel operators and their applications, Literature [8] [9] is the work of composition operators, while Zhu's monograph [10] is the classic work of operator theory on function space, and literature [13] comprehensively introduces operator theory, which can be used as an introduction to operator theory.
\par  Composition operator~$C_{\varphi}$~ on Hilbert spaces never turns out. The related work of ~$C_{\varphi}$~ is a very large, interesting, and profound problem. Naturally, we can associate their algebraic properties, invertibility compactness and even eigenvalue spectrum, etc. In general, we are very glad to develop the knowledge about composition operators~$C_{\varphi}$~ on some Hilbert space,they features many unique properties.An extended eigenvalue~$\lambda$~ is an algebric element of a bounded linear operator A which satisfies there exsits nonzero operator X such that ~$$AX=\lambda XA$$~ on a Hilbert space, we discuss the extended eigenvalues of composite operators. For a bounded linear operator, If the nonzero operator~$X $~and the nonzero complex scalar~$ \lambda $~exist on the Hilbert space, which makes~$AX=\lambda XA $~hold, then we will call~$ \lambda $~the extended eigenvalue of this operator, and~$X$~the extended eigenvalue operator. In , several years ago,Miguel is the most important role in the feild,he tramendously extended the analysis in extended eigenvalues of composition operators,and he totally illuminated the extended eigenvalues for composition operators induced on the Hardy space by linear fractional transformations of the unit disk.He made the key to lead to intriguing questions on the interface between complex analysis and operator theory. Here are the conclusions.Our trip start by him and finally get some useful answers on Bergman spaces and Fock spaces. The extended eigenvalues of composition operators on Fock spaces are under our control. While most work has been pushed on Bergman spaces, but unfortunately others have not been totally answered yet.
\subsection {Extended eigenvalues and extended eigenoperators}

\par Premising that the infinite dimensional complex separable Hilbert space~$H$~. As mentioned above, the complex scalar~$\lambda$~is called the extended eigenvalue of the bounded linear operator~$A$~.If there is a non-zero operator~$X$~, so that~$AX=\lambda{XA}$~, then~$ X $~is called the extended eigenoperator of operator~$A$~. The set consists of all extended eigenvalues of operator~$A$~is called the extended spectrum of~$A$~, what we called ~$Ext(A)$~. In addition, the symbol~$\varepsilon{xt}(\lambda,A)$~represents the set of all extended eigenoperators under a given extended eigenvalue~$\lambda$~. If zero operator is an extended eigenoperator, the set of extended eigenoperators will become a weakly closed linear manifold. These previous studies are derived from the works of Scott brown[18] and Kim, Moore and pearcy[19]. Victor lomonosov[16] proved that if an operator~$A$~has a non-zero compact extended eigenoperator, then the operator~$A$~has an extraordinary closed hypervariant subspace. The systematic study of extended eigenvalues and extended characteristic operators of classical operators began with Biswas, Lambert [14] worked on the Volterra operator and will continue to work on Lambert [15]. He also gave the sufficient conditions for the operator to have a nontrivial closed highly variable subspace and who calculated the extended spectrum of scalar perturbation of the Volterra operator. At present, this research field has become a very active research field, many contributions have been made in the invariant subspace [17,18] and the calculation of extended spectra of some special class operators [14,15,16]
\subsection {Composition operators on Hardy, Fock, Bergman Spaces}

The composition operator what we will discuss is based on Hardy space~$H^{2}{(\mathbb{D})}$~, Bergman space~$L^{2}_a$~and Fock space~$F^{2}_{\alpha}$~. The analytic self mapping of the unit disk~$\varphi$~expresses the linear operator~$C_{\varphi}$~on the space of analytic functions on the unit disk or complex plane by the following expression:~$$C_{\varphi}f=f\circ\varphi, f:\mathbb{D}\rightarrow\mathbb{C} $$~ is analytical.
Littlewood's membership principle [27] ensures that the composition operator is a bounded linear operator on Hardy space~$H^{2}{(\mathbb{D})}$~. The main reason is that~$H^{2}{(\mathbb{D})}$~is the most natural Hilbert space in analytic functions. From the Fatou lemma, we can see that the mapping of function~$f$~on each~$H^{2}{(\mathbb{D})}$~to its radial limit f is almost everywhere represented by the expression~$f(e^{i\theta})=\lim\limits_{r\to1}f(re^{i\theta})$~decision is an isometric mapping from~$H^{2}(\mathbb{D})$~to~$L^{2}(\mathbb{T})$~. Therefore, scholars have developed a complete set of theories around the functional and operator theoretical properties of composition operators on Hardy spaces. Joel Shapiro's book [27] is a very profound paper on this topic. The theoretical system of ~$A$~is based on Hardy space. According to some previous conclusions, the corresponding relationship between the boundedness of composition operators and their index functions in Bergman space~$L^{2}_a$~ is never different from Hardy space[25]. Zhu told us discussing the boundedness of composition operators in Bergman space are not supposed to. However, there are some differences about Fock space. The previous conclusions show that the boundedness of composition operators in Bergman space~$L^{2}_a$~ is similar to that in Hardy space If and only if~$$\varphi(z)=az+b, a,b\in\mathbb{C},|a|\leq1$$~is a bounded linear operator on Fock space[20]. Fock space has more particularity than other spaces, especially goes to the boundness.
\subsection {classification of linear fractional transformation}
 Recall the linear fractional mapping~$$f(z)=\frac{az+b}{cz+d}$$~defined on the Riemannian sphere~$\widehat{\mathbb{C}}=\mathbb{C}\bigcup[{+\infty}]$~. The condition~$ad-dc=0$~ensures that~$f$~is not a constant function. Usually we need to consider the infinite point. As we all know,~$f$~is an automorphism of the Riemannian sphere, and from a metric point of view, ~$f$~is also an isometry. A necessary and sufficient condition for the linear fractional mapping~$f$~to be self mapping is that~$|bd-ac|\leq |b|^{2}-|d|^{2}$~. 
\par This simple case of par rarely appears in the literature. For example, it can be found in the works of Maríaj.martínor [24]in Contreras, DíAZ Madrid, martíN and vukotić.The latter reference provides another condition, which is more complex to some extent, but has an easier proof.
The general discussion of LFT is very complex.
We usually classify LFT according to its fixed point situation and geometric properties. Refer to Joel shapiro for classification[20]. In Miguel'case[27], there are eight important LFTs: elliptic automorphism(EA), parabolic automorphism(PA), parabolic non automorphism(PNA), hyperbolic automorphism(HA), class I, class II or class III hyperbolic non automorphism,(HNA I)(HNA II)(HNAIII)and loxodromotic mapping(LOX). After appropriate linear fractional mapping, each LFT can be represented in the standard form collected in Figure1. 
\par In Shapiro's classification, there is no HNAIII class of hyperbolic non automorphic LFT without fixed point on~$\mathbb{T} $~, but there is a simple solution to this problem, because these mappings have the same standard form as isometric mapping.

\subsection {quasi-differential operator and quasi-integral operator}
\par Quasi-differential operators and quasi-integral operators do not pay much attention to Hardy and Bergman spaces~$L^{2}_a$~ but  play an important role about Fock spaces. Our research progress is also related to them on ~$F^{2}_{\alpha}$~,  define two operators~$D$~and~$X$~, where~$$Xf(x)=xF(x), DF(x)=f'(x).$$~ The style of quasi-differential operator and quasi-integral operator is very simple. Generally, it is hard to discuss their properties on Bergman Spaces~$L^{2}_a$~ and Hardy spaces. They are unbounded dense operators. But on the subspace of the holomorphic function space Fock space, the discussion of unbounded operators will be very natural.Toeplitz Operators in Fock spaces can be completely characterized by quasi-differential operators and quasi-integral operators[].
\section{preliminary results}
\par In this section, we  establish several lemmas which will run through this paper.
\begin{Lemma}\label {lemma1}Operator~$A\in{B(H)} $~, then it is injective if and only if~$0\notin{Ext(A)} .$~
\end{Lemma}
This can be proved by the open mapping theorem and the analytic continuation principle.
\begin{Lemma}
\label {lemma2} We assume that~$\varphi$~is a non constant analytic self map of the open unit disk, then~$C_\varphi$~is a injective linear operator.
\end{Lemma}
If we combine lemma 3.1 and lemma 3.2, directly conclude that 0 is never the extended eigenvalue of the composition operator.
\begin{Lemma}
\label {lemma3} If~$\varphi$~is a non constant analytic self mapping of the open unit disk, then~$0\notin{Ext(C_\varphi)}$~ \end{Lemma}
Next, we give two basic results of extended spectrum under adjoint and scalar multiplication
\begin{Lemma}
\label {lemma4} If~$A\in{B(H)}$~, then~$Ext(A^*)/{0}={1/\overline\lambda}$~ \end{Lemma}
\begin{Lemma}
\label {lemma5} If~$A\in{B(H)}$~and~$\alpha\in\mathbb{C}/{0} $~, then~$Ext(\alpha{A})=Ext(A)$~
\end{Lemma}
The following result is a necessary condition for the extended eigenvalue of a injective operator with a set of eigenvectors. We will use~$\sigma_ {p}(A)$~denotes the point spectrum of the operator~$A\in{B(H)} $~, that is, the set of all eigenvalues of~$A$~.
\begin{Lemma}
\label {lemma6} Let~$A\in {B(H)} $~be a injective operator with non-empty point spectrum. If the set of eigenvectors of~$A$~can form a dense subset of Hilbert space, then~$Ext(A)\subseteq {\{ \alpha/\beta:\alpha, \beta\in\sigma_{p}(A)\}} $~
\end{Lemma}
{\bf Proof.}Let~$F$~be a total subset of~$H$~consisting of the eigenvectors of~$A$~. Let~$\lambda\in Ext(A)$~and~$X\in\varepsilon{xt}{(\lambda, a)} $~. Secondly, there is~$f\in F$~such that~$Xf\neq0$~, and there is~$\beta\in\sigma_ {p}(A)$~make~$Af=\beta{f}$~
From this, we get~$$AXf=\lambda XAf=\lambda\beta{Xf}. $$~ Since~$Xf\neq0$~, we know ~$$\alpha:=\lambda\beta\in\sigma_{p}(A).$$~ Since ~$A$~is injective and we have~$\beta\neq0$~, ~$$\lambda=\alpha/\beta,$$~ just what we want.\\
The multiplication operator is easily associated with the extended eigenvalue of the composition operator. Every bounded analytic function~$b\in{H^{\infty}(\mathbb{D})} $~derives ~$M_{b} $~Multiplication operator~$M_{b}$~defined by ~$$M_{b}f=bf,$$~ ~$M_{b}\in B(H^{2}(\mathbb{D}))$~. Easy to check~$M_{b} $~is indeed a bounded linear operator.
$\hfill\square$
\begin{Lemma}\label{lemma7} set~$C_\varphi$~is a composition operator induced by~$\varphi$~. Let~$\lambda\in\sigma_{p}(C_\varphi)$~,~$b$~are functions in~$H^{\infty}(\mathbb{D})$~, so that~$C_\varphi{b}=\lambda b$~. And yes,~$\lambda\in Ext(C_{\varphi})$~and~$C_{\varphi}M_B=\lambda M_{b} C_{\varphi}$~
\end{Lemma}
{\bf Proof.}Note that for each~$f\in H^ {2}(\mathbb{D})$~, we have~$$C_{\varphi}M_{b}f=C_{\varphi}(b\circ{f})=(b\circ\varphi)(f\circ\varphi)=(C_{\varphi}b)(C_\varphi{f})=\lambda{b}C_{\varphi}f=\lambda{M_{b}}C_{\varphi}f $$~
, which makes~$$C_{\varphi}M_ B=\lambda M_ {b} C_ \varphi.$$~\\
Easy to see this lemma has everything to do also about Bergman space~$L^{2}_a$~.But we can't apply it in Fock space since the multiplication operator~$M_{b}$~ cannot be a bounded linear operator on Fock space.
$\hfill\square$\\
The last two results of the extended spectrum of the direct sum of two operators in this section will help to simplify the problem of the second kind of hyperbolic non automorphic composition operators(HNA II)to the problem of the first kind(HNA I).
\begin{Lemma}
\label {lemma8} consider a direct sum decomposition~$H=H_{1} \oplus{H_2}$~, set~$A_{1} \in B(H_1), A_{2} \in B(H_2)$~. then~$Ext(A_1)\cup{Ext(A_2)}\subseteq{(A_{1}\oplus{A_2})} $~ \end{Lemma}
\begin{Lemma}
\label {lemma9} consider a direct sum decomposition~$H=H_{1} \oplus{H_2}$~, set~$I_1$~is~$H $~. Let~$A$~be a injective operator in~$B(H_2)$~,
~$Ext(I_{1}\oplus{A})=\sigma_{p}(A^*)\bigcup{{\overline{\sigma_{p}(A)}} ^ {-1}}\bigcup Ext(A^*)$~ .\end{Lemma}
\par Lemmas 3.1 to 3.3 emphasize that 0 does not become an element in the extended spectrum of composite operators, which may seem minor, but in reality, their existence is necessary because once 0 elements are in the extended spectrum, it can cause significant problems in Lemma 2.6 and directly block our research. Lemmas 2.4 and 2.5 actually supplement the influence of scalars on the extended spectrum of composite operators. As a special algebraic property, the extended spectrum is similar to other algebraic structures in that they are star corresponding and scalar invariant. They can help reduce the details in the proof and simplify the proof process. Lemma 2.6 is one of the core ideas that runs through the entire article. Broadly speaking, it is relatively easy for us to find the spectrum of an operator. And its extended spectrum is difficult for us to estimate or calculate, but Lemma 2.6 can connect the spectrum with the extended spectrum. Furthermore, if the target operator satisfies the condition of 2.6, we can directly obtain an "upper bound", or boundary, of its extended spectrum.
 \par  If 2.6 is about finding the upper limit, then the following theorem is a means of finding the lower limit. Theorem 2.7 tells us a special conclusion when the characteristic function is a bounded analytic function. Although 2.7 does not have universality, this property is extremely effective for composite operators in Bergman spaces. Assuming that all the operators we explore satisfy the condition of Lemma 2.6 (in fact, all the operators discussed in this article are like this, even those operators for which we have not reached specific conclusions), a natural idea arises, if we can find a way to prove that all the elements in the "upper bound" satisfy the property of extended eigenvalues, or directly use algebraic and topological properties to prove that the extended spectrum is a closed set, then all our problems will be solved easily. Therefore, it seems that we have found the Eureka to solve the problem, and the rest is to complete the entire process or verify various feasibility.
\par Firstly and most importantly, we want to use the method of analogizing Hardy spaces to handle the situation of Fock and Bergman spaces~$L^{2}_a$~, which is also the most normal and possible approach. After in-depth exploration and research, we can summarize the key elements of calculating the extended spectrum as follows: finding the point spectrum and eigenvalues of the operator, verifying the conditions of Lemma 2.6, calculating the upper bound, finding a way to calculate the lower bound, and calculating the lower bound. The key point of the validation of Lemma 2.6 is whether it constitutes a total subset. The validation of some operators is relatively simple, and we can solve the slightly complex problem by introducing the concept of rich point spectrum. Relatively speaking, the determination of the upper bound is not complicated or cumbersome, and the difficulty lies in the lower bound. We fantasized about completely characterizing the lower bound through lemmas 2.7-2.9, but it went against our expectations and not all of my research was smooth sailing. In fact, if we consider the case of Fock space, it is easy to find that according to the Liouville theorem, a function must be a constant value function if it is a bounded linear mapping in the entire complex plane. Therefore, our discussion becomes meaningless. So we need to find ways to improve the proof approach in the Fock space.\par In the first part, we have discussed that there are only a few fractional linear mappings on Fock spaces that can induce bounded composite operators. We first introduce the relevant properties obtained by the composite operator of elliptic self mapping. As is well known, weight is usually a factor worth considering in Fock space. On Fock space and Fock type weighted bounded analytic function space, if the linear fractional transformation satisfies~$\varphi(z)=w{z}$~where~$|w|=1$~, then we consider the weighted multiplication operator ~$K_{b}$~, and have the following conclusion:
\begin{Lemma}\label{lemma10} set~$C_\varphi$~is a composition operator induced by~$\varphi$~on Fock space. Let~$\lambda\in\sigma_{p}(C_\varphi)$~,~$b$~are functions in~$F^{\infty}_{\alpha}$~, then~$C_\varphi{b}=\lambda b$~. And~$\lambda\in Ext(C_{\varphi})$~and~$C_{\varphi}K_b=\lambda K_{b} C_{\varphi}.$~
\end{Lemma}
{\bf Proof.}The proof of lemma 3.10 is similar to 2.7, the only point we need to change is that we must use the weighted multiplication operator ~$K_{b}$~ instead of the multiplication operator~$M_{b}$~, which is defined as ~$$K_{b}f=b(z)f(z)e^{-\alpha{|z|^2}}$$~We note that for each~$f\in F^{2}_{\alpha}$~, we have~$$C_{\varphi}K_{b}f=C_{\varphi}(b\circ{f}{\circ}e^{-\alpha{|z|^2}})=(b\circ\varphi)(f\circ\varphi)e^{-\alpha{{w}|z|^2}}$$~~$$=(C_{\varphi}b)(C_\varphi{f})e^{-\alpha{|z|^2}}=\lambda{b}e^{-\alpha{|z|^2}}C_{\varphi}f=\lambda{K_{b}}C_{\varphi}f$$~
It makes~$$C_{\varphi}K_ b=\lambda K_ {b} C_ \varphi.$$~
$\hfill\square$
\par Trying the similar ways for other type composition operators on Fock space, we fail to achieve our wide research on Fock spaces, but we find that the quasi-differential operator~$D $~and the quasi-integral operator~$X$~play a key role. On ~$F^{2}_{\alpha}$~, we define~$D $~and~$X $~as follows,~$$Xf (x)=xf (x), Df (x)=\frac {1} {2\alpha {i} }f '(x) $$~. Using the quasi-differential operator~$D $~ and the integral operator~$X $~. Unfortunately, both quasi-differential operators and quasi-integral operators are unbounded dense  operators and not every element in Fock space can be differential infinitely. Thanks to the atomic decomposition theorem on Fock space,  the normalized reproducing kernel function~$k_ {z} $~, the linear combination of ~$e ^ { \alpha {w}\overline{z}} $~is dense in~$F^ {2}_{\alpha} $~, both quasi-differential operators and quasi-integral operators can operate on the dense set of the normalized reproducing kernel function~$k_ {z} $~in~$F^{2}_{\alpha}$~.More discussions will be in discussed next chapter.

\section{Extended eigenvalues on Fock spaces}
Recall that if ~$\varphi$~ is an elliptic automorphism, then~$ \varphi$~has two fixed points, one in ~$\mathbb{D}$~ and another in ~$\mathbb{C}\setminus\mathbb{D}$~. Given this fixed point configuration, upon conjugation by disk automorphisms, we can suppose that the
fixed points are 0,~${\infty}$~, which leads to the standard form ~$\varphi(z)=w{z},|w|=1$~. Let us observe that the conjugation by disk automorphisms induces a similarity at the operator level, and the extended spectra are invariant under similarities. Notice that ~$w_n$~ is an eigenvalue of ~$C_{\varphi}$~, and a corresponding eigenfunction is the monomial ~$e_{n}(z)=z^n $~, for every ~$n\in \mathbb{N}_0.$~

\subsection{ Elliptic automophism}
\begin{Theorem}\label{Theorem1}
{If~$\varphi(z)=w{z}$~where~$|w|=1$~,~$Ext(C_\varphi)=\{ {w}^n:n\in\mathbb{Z}\}$~}.\end{Theorem}
It is completely consistent with the proof on Hardy space, and the extended feature operator is also the same. What we wll use is still a left shift operator on Hardy space, but not same on Bergman space. In this case, it is still a bounded linear operator.In Fock type spaces,we often call shift operators as Mobius invirant maps.Lemma 3.10 and Lemma 2.6 are the key to proving theorem 4.1. We first use lemma 2.6 to determine the upper bound of the extended spectrum, and then use lemma 3.10 to obtain some elements in the extended spectrum, which are still included in this upper bound.
From culculations we have~$[e_{n}:n\in\mathbb{N}_{0}]$~are the eigenvalues set of~$C_{\varphi}$~,all we know that basic vectors in Fock have the types:~$$e_{n}(z)=\sqrt{\frac{{\alpha}^{n}}{n!}}z^n$$~,then~$[e_{n}:n\in\mathbb{N}_{0}]$~is a total family,directly use lemma 2.6 we have ~$Ext(C_\varphi)\subseteq[{w}^n:n\in\mathbb{Z}]$~.
\begin{Theorem}\label{Theorem2} {If~$\varphi(z)=w{z},|w|=1$~,~$Ext(C_\varphi)\subseteq\{ {w}^n:n\in\mathbb{N}_{0}\} $~}
\end{Theorem}
Theorem 4.2 is easily gained by lemma 3.10. Now we need to discuss~$[{w}^{-n}:n\in\mathbb{N}_{0}]$~.
\begin{Theorem}\label{Theorem3} {If~$\varphi(z)=w{z},|w|=1$~,~$Ext(C_\varphi)\subseteq\{{w}^{-n}:n\in\mathbb{N}_{0}\}$~}
\end{Theorem}
{\bf Proof.}We define S as ~$$M^{*}_{e_k}M_{e^{-\alpha|z|^2}}f(z)=Sf(z)$$~,~$$C_{\varphi}{Se_{n}(z)}=\sqrt{\frac{{\alpha}^n}{n!}}C_{\varphi}{Sz^n}=\sqrt{\frac{{\alpha}^n}{n!}}C_{\varphi}{z^{n-k}e^{-\alpha|z|^2}}=\sqrt{\frac{{\alpha}^n}{n!}}{w}^{n-k}z^{n-k}e^{-\alpha|z|^2}$$~~$$=\frac{\sqrt{\frac{{\alpha}^n}{n!}}}{\sqrt{\frac{{\alpha}^{n-k+1}}{(n-k+1)!}}}{w}^{n-k}e_{n-k}(z)e^{-\alpha|z|^2}={w}^{-k}SC_\varphi{e_{n}(z)}$$~for every n and k.~$M^{*}_{e_k}$~ is still a bounded linear operator on Fock space.And from the well-known truth that polynomials are dense in Fock space we get theorem 4.3.\\
Combine theorem 4.2 and 4.3 we can easily have the conclusion ~$Ext(C_\varphi)=\{{w}^n:n\in\mathbb{Z}\}$~,so that we finish part one.
$\hfill\square$
\subsection{Other types}
Following we have completed everything in the type of elliptic automophism.In other situation, what we called other types are functions like ~$\varphi(z)=az+b, a,b\in\mathbb{C},|a|<1 $~ or  ~$|a|=1$~  with  ~$b=0$~.
\par Naturally we have did everything in this field,
\begin{Theorem}\label{Theorem4}
{If~$\varphi(z)=w{z}+b,|w|\leq1,b\in\mathbb{C}$~,~$Ext(C_\varphi)=\{{w}^n:n\in\mathbb{Z}\}$~}\end{Theorem}
Theorem 4.4 contained all types of composition operators on Fock space.At first we naturally consider the same way we have used before, it has been proved that  composition operator ~$C_{\varphi}$~ has a countable point spectrum, namely ~$\sigma_{p}(C_{\varphi}) =\{{\varphi} '(c)^n : n\in\mathbb{N_{0}}\}$~. Further, each eigenvalue is of multiplicity one and a corresponding eigenfunction is given by~$ h(z) = z-\frac{b}{w}$~
\begin{Theorem}\label{Theorem5}
{If~$\varphi(z)=w{z}+b,a,b\in\mathbb{C},0<|w|<1 $~,then~$Ext(C_\varphi)\subseteq\{{w}^n:n\in\mathbb{Z}\}$~}\end{Theorem}
Continue our discussion under the eigenfuctions and spectrums, the given eigenfunction ~$ h(z) = z-\frac{b}{w}$~ can be spanned to be a  dense manifold in Fock space.With lemma3.5 and ~$\sigma_{p}(C_{\varphi}) =\{{\varphi} '(c)^n : n\in\mathbb{{N}_{0}}\}$~ we get theorem 4.5.
\begin{Theorem}\label{Theorem6}
{If~$\varphi(z)=w{z}+b,a,b\in\mathbb{C},0<|w|<1 $~,then~$\{{w}^n:n\in\mathbb{Z}\}\subseteq Ext(C_\varphi)$~}\end{Theorem}
Back to the defination.Such as lemma 3.10 or 2.7,we try to bulid ~$$K_{b}f=b(z)f(z)e^{-\alpha{|z|^2}}$$~ samely, but it cannot be possible since calculations that ~$C_{\varphi}K_{b}f=\lambda{K_{b}}C_{\varphi}f$~will only hold if~$\varphi$~elliptic maps.The reason is that ~$e^{-\alpha{|z|^2}}$~is invirant only under the elliptic map,then we get~$\{{w}^{-n}:n\in\mathbb{N}\}\subseteq Ext(C_\varphi)$~.So the key to solve our problem absolutely is to find out what is invirant under the fractional transform, we are very glad to see the invirant elements with some stable features.\\

In this part we need the new theory to continue our process, now the important point that we discuss comes to the quasi-differential operator~$D $~ and the quasi-integral operator~$X $~.~$$Xf (x)=xf (x), Df (x)=\frac {1} {2\alpha {i} }f '(x) $$~,then we directly have ~$$C_{\varphi}Df(z)=C_{\varphi}(\frac {1} {2\alpha {i} }f'(z))=\frac {1} {2\alpha {i} }{w}^{-1}f'({\varphi}(z))={w}^{-1}DC_{\varphi}f$$~and mutiplitively~$$C_{\varphi}D^{n}f(z)=C_{\varphi}((\frac {1} {2\alpha {i} })^{n}f^{(n)}(z))=(\frac {1} {2\alpha {i} })^{n}{w}^{-n}f^{(n)}({\varphi}(z))={w}^{-n}D^{n}C_{\varphi}f$$~.
\par The key road we need to go with is that some elements on Fock space aren't differential. Theorem 2.11 in ~$Fock spaces$~ by Zhu told us functions like[]are dense on Fock space and they are infinitesimal.The using of the lemma from Zhu makes the proof comes to be completed.Finally we finish theorem 4.4 by the combination of theorem4.5 and 4.6.
\par The following work illuminate the extended spectrum of composition operators on Fock type space,we can directly conclude that the  extended spectrum of composition operators on Fock type space are tremendously consist of the scalar constant of the linear fractional transform~$\varphi$~.
\par Fock spaces have many different features with Hardy spaces and Bergman spaces,theory about ~$H^{\infty}{(\mathbb{D})}$~ do many things with Hardy spaces and Bergman spaces. Garnett exactly illuminate them. However, we can easily notice that ~$H^{\infty}{(\mathbb{D})}$~ have less thing to do with Fock type spaces. To our luck, we deal the problem smoothly on Fock type spaces provided that composition operators are easily to be handled. Now, the following discussion will be back to traditional analysis about Bergman spaces which don't mean that there are nothing left creatively. We still own some small thought about situation on Bergman spaces.

\section{Extended eigenvalues on Bergman spaces~$L^{2}_a$~}

\subsection{ Elliptic automorphism}

 \begin{Proposition}\label{proposition1}
{If~$\varphi(z)=w{z},|w|=1$~,~$\{{w}^n:n\in\mathbb{N}\}\subseteq Ext(C_\varphi)\subseteq\{{w}^n:n\in\mathbb{Z}\}$~}
\end{Proposition}
{\bf Proof.}It is completely consistent with the proof on Hardy space, and the extended feature operator is also the same. It is still a left shift operator on Hardy space, but not a left shift operator on Bergman space~$L^{2}_a$~. In this case, it is still a bounded linear operator so that it would be useful probably.  Lemma 2.7 is the key point to understand the spectrum in Hardy space, and it also holds on Bergman spaces .
So we have the thought. Notice that ~$w$~ is an eigenvalue of composition operator, and a corresponding eigenfunction is the monomial ~$e_{n}(z) = z^n $~, for every ~$n\in \mathbb{N}_0$~.	On the other hand, the set ~${e_n : n\in \mathbb{N}_0} $~of eigenfunctions for ~$Ext(C_\varphi)$~ is a total subset of ~$L^{2}_a$~.Then a straight use of lemma2.6 and lemma2.7 we get theorem1. 
$\hfill\square$\\
Observe the theorem, what we remind is the relationship that the posibility of~$\{{w}^{-n}:n\in\mathbb{N}\} \subseteq Ext(C_\varphi)$~.
 \begin{Theorem}\label{Theorem7}
{If~$\varphi(z)=w{z},|w|=1$~,~$\{{w}^{-n}:n\in\mathbb{N}\} \subseteq Ext(C_\varphi)\subseteq\{{w}^n:n\in\mathbb{Z}\}$~}
\end{Theorem}
{\bf Proof.}Consider the shift operator X for a fixed positive integal k as:~$$Xe_n  = e_{n-k}$$~  for ~$n\neq k$~,  otherwise~$Xe_n=0$~ . Zhu's []told us ~$X$~ is a bounded linear operator on ~$L^{2}_a$~ for every fixed k. In this case, if ~$n\neq k$~,~$$C_{\varphi}X_k{z^n}=C_{\varphi}z^{n-k}={w}^{n-k}z^{n-k}={w}^{n-k}X_k{z^n}={w}^{-k}X_{k}C_{\varphi}{z^n}$$~ In other time we have ~$$C_{\varphi}X_k{z^n}=0={w}^{-k}X_{k}C_{\varphi}{z^n}(k>n)$$~or ~$$C_{\varphi}X_k{1}=0={w}^{-k}X_{k}C_{\varphi}{1}, n=0$$~ Anyway we prove that ~$$C_{\varphi}X_k{z^n}={w}^{-k}X_{k}C_{\varphi}{z^n}$$~ We know that ~$[1,z,z^2,z^3...z^n...]$~is a dense subset of Bergman space. So for every function f in ~$L^{2}_a$~,~$there exists p as a polynomial, s.t.\Vert{f-p}\Vert \rightarrow 0$~, and ~$$\Vert{C_{\varphi}X_k{f}-{w}^{-k}X_{k}C_{\varphi}f}\Vert\leq\Vert{C_{\varphi}X_k{f-p}}\Vert +\Vert{{w}^{-k}X_{k}C_{\varphi}{f-p}} \Vert\rightarrow 0$$~ And with the help of density of polynomials in ~$L^{2}_a$~ we get theorem 5.2.
$\hfill\square$
\begin{Theorem}\label{Theorem8}
{If~$\varphi(z)=w{z},|w|=1$~,~$Ext(C_\varphi)=[{w}^n:n\in\mathbb{Z}]$~}\end{Theorem}
  
\par Combine 5.1 and 5.2 we have the result.In this time the extended eigenoperator is the left shift operator X.

\subsection{loxodromic map or a hyperbolic nonautomorphism of the third kind}
\par Last year, Bensaid calculated the extended spectrum of equiangular or HNAIII on a weighted Hardy space, and the ~$L^{2}_a$~ happens to be a special case of the weighted Hardy space.
Recall that if ~${\varphi}$~is a loxodromic map or a hyperbolic nonautomorphism of the third kind, ~${\varphi}$~can be assumed to have one fixed point at infinity, while the other one, say c, belongs to ~$\mathbb{D}$~. In this case ~${\varphi}$~has the standard form that ~$\varphi=a(z-c)+c$~, the relative eigenvalues are ~$\{\varphi'(c):n\in\mathbb{N}\}$~, and the eigenfunctions are ~$\{\sigma^{n}(z):\sigma(z)=z-c$~ where ~$n\in\mathbb{N}\}$~.The hyperbolic nonautomorphism of the third kind corresponds to the case ~$a>0$~. The composition operator ~$C_{\varphi}$~has a countable point spectrum this time, namely ~$\sigma_{p}(C_\varphi)=\{{{\varphi'(c)}^n : n\in\mathbb{N}_{0}}\}$~. Furthermore, each eigenvalue is of multiplicity one and a corresponding eigenfunction is given by ~${\sigma}^n=(a(z-c)+c)^n$~
. \par So~$$C_{\varphi}{\sigma}^n={\varphi(c)}^{n}{\sigma}^n$$~  The span of the eigenfunctions~$\{{{\sigma}^n : n\in\mathbb{N}_{0}}\}$~ is a dense linear manifold in ~$H^{2}(\mathbb{D})$~ and~$L^{2}_a$~.  All eigenfunctions are bounded on Bergman space and Hardy space.This means that we can get the next theorem by lemma 2.7 and 2.6.
\begin{Theorem}\label{Theorem9}
{If~$\varphi(z)=a({z}-c)+c,c\neq0,|a|+|1-a||c|\leq1~$ or$~ |a|<1,c=0$~, ~$\varphi$~is loxdromotic or HNAIII,then~$\{{{\varphi}'(c)}^n:n\in\mathbb{Z}\}\subseteq Ext(C_\varphi)\subseteq\{{{\varphi}'(c)}^n:n\in\mathbb{Z}\}$~}\end{Theorem}
\par The following theory is an easy application of lemmas in chapter 3. What we are supposed to do is whether~$$\{{{\varphi}'(c)}^{-n}:n\in\mathbb{Z}\}\subseteq Ext(C_\varphi).$$~
\begin{Theorem}\label{Theorem10}
{If~$\varphi(z)=a({z}-c)+c, c\neq0,|a|+|1-a||c|\leq1~$ or$~ |a|<1, c=0$~, ~$\varphi$~is loxodromotic or HNAIII, then~$Ext(C_\varphi)=\{{{\varphi}'(c)}^n:n\in\mathbb{Z}\}$~}\end{Theorem}
Consider the Toeplitz operator~$T_{{\sigma}^{-1}}$~,for every function ~$f\in L^{2}_a$~, ~$$T_{{\sigma}^{-1}}f=P({\sigma}^{-1}f)$$~.~$P$~is the analytic mapping from ~$l_2$~ to~$L^{2}_a$~.In this case, if ~$n>1$~, ~$$C_{\varphi}T_{{\sigma}^{-1}}{{\sigma}^n}=C_{\varphi}{\sigma}^{n-1}={\varphi('c)}^{n-1}{\sigma}^{n-1}={\varphi'(c)}^{n-1}T_{{\sigma}^{-1}}{{\sigma}^n}={\varphi'(c)}^{-1}T_{{\sigma}^{-1}}C_{\varphi}{{\sigma}^n}.$$~  If ~$k<n$~ and ~$k\in\mathbb{N}$~,~$$C_{\varphi}T_{{\sigma}^{-k}}{{\sigma}^n}=C_{\varphi}{\sigma}^{n-k}={\varphi('c)}^{n-k}{\sigma}^{n-k}={\varphi'(c)}^{n-k}T_{{\sigma}^{-k}}{{\sigma}^n}={\varphi'(c)}^{-k}T_{{\sigma}^{-k}}C_{\varphi}{{\sigma}^n}.$$~ And if ~$n=1$~, ~$k$~ is of empty. Otherwise, the series ~$\{{\sigma}^{n}:n\in\mathbb{N}\}$~is a dense subset of ~$L^{2}_a$~. Like the way we used in last part we directly get for every ~$f\in L^{2}_a, k>0$~, ~$$C_{\varphi}T_{{\sigma}^{-k}}f={\varphi'(c)}^{-k}T_{{\sigma}^{-k}}C_{\varphi}f$$~ Which means ~$$C_{\varphi}T_{{\sigma}^{-k}}={\varphi'(c)}^{-k}T_{{\sigma}^{-k}}C_{\varphi},$$~and ~$$\{{{\varphi}'(c)}^{-n}:n\in\mathbb{Z}\}\subseteq Ext(C_\varphi)$$~. Combine it with theorem 5.4 we finally get ~$Ext(C_\varphi)=\{{{\varphi}'(c)}^n:n\in\mathbb{Z}\}$~.\\

\subsection{hyperbolic automorphism}
\par If ~${\varphi}$~is a hyperbolic automorphism of the unit disk, its fixed points lie on the boundary of the unit circle, and upon conjugation by disk automorphisms, we can suppose that the fixed point of ~${\varphi}$~are  1 and -1.
\par Moreover, we can suppose that ~${\varphi}$~has the following standard form.~$$\varphi(z)=\frac{z+r}{1+rz},0<r<1.$$~ In what follows we will take advantage of some spectral information about~$C_{\varphi}$~that can be found in [], every point in the ring ~$$G:=\{\alpha\in\mathbb{C}:R^{-1/2}<|\alpha|<R^{1/2}\}$$~that ~$$R=\frac{1+r}{1-r}$$~can be regarded as an eigenvalue for ~$C_\varphi$~, 
the correspanding eigenfunction is ~$$e_{w}(z)=({\frac{1+z}{1-z}})^{w}.$$~This time we have the formula ~$$C_{\varphi}e_{w}(z)=R^{w}e_{w}(z).$$~
The notion of an operator with rich point spectrum has been introduced recently as a way to determine the extended eigenvalues for Cesàro operators [14] and bilateral weighted shifts [13,14].\par 
An operator ~$A\in B(H)$~is said to have a rich point spectrum provided that ~$int\sigma_{p}(A)\neq $~ and for every open disk ~$$D\subseteq\sigma_{p}(A)$$~, the family of the corresponding eigenvectors.In this situation~$$\bigcup ker(A-z)	$$~is a total subset of ~$H$~. The following result []is a suﬃcient condition for an operator ~$A $~to have rich point spectrum.
\begin{Lemma}\label{lemma14}Let  ~$ A\in B(H)$~and assume there are an open connected set ~$ w\subseteq\mathbb{C}$~, an analytic mapping ~$h : w\rightarrow H$~ and a nonconstant analytic function ~$\gamma: w\rightarrow\mathbb{C} $~such that:\\
(i)	~$h(w)\in ker[A-\gamma(w)] /{0}$~ for all~$ w\in w$~,\\
(ii)	~${h(w): w\in w} $~is a total subset of ~$H$~, \\
(iii)	~$\sigma_p(A)\subseteq clos\gamma(w)$~.\\
Then A has rich point spectrum.
\end{Lemma}
\par It turns out that lemma applies to the composition operators under consideration. There are some deep lemmas we'd use step by step.
\begin{Lemma}\label{lemma23}The correspanding eigenfunction is ~$$e_{w}(z)=({\frac{1+z}{1-z}})^{w}$$~ is invertible on~$L^{2}_a$~, moreover it is outer and cyclic, the set ~$\{e_{w}:w\in K\}$~ is a total set.
\end{Lemma}
{\bf Proof.}:Observing ~$e_{w}(z)$~ we find that ~$e_{w}(z), {e_{w}(z)}^{-1} $~are both in ~$H^{\infty}{(\mathbb{D})}$~, absolutely in Bergman space. By Zhu's lemma, an element is invertible means that it is cyclic, cyclic vectors are same as outer vectors on Bergman space.
$\hfill\square$
\begin{Lemma}\label{lemma18}If ~${\varphi}$~is a hyperbolic automorphism of the unit disk,
then~$C_{\varphi}$~has rich point spectrum.
\end{Lemma}
{\bf Proof.}Let ~$A=C_{\varphi}$~and ~$h(w):=e_{w} $~and~$ \gamma=R^{w}$~.First, ~$$C_{\varphi}e_{w}(z)=R^{w}e_{w}(z)$$~ means condition (i) holds.The proof of (ii) has been introduced before.Finally, it is clearly that~$w$~ is an open connected set and ~$ \gamma$~ is analytic such that ~$\gamma=G=\sigma_{p}(C_{\varphi})$~,which solve condition(iii).
$\hfill\square$
\begin{Lemma}\label{lemma20}Let  ~$ A\in B(H)$~ and A has rich point spectrum, there are constants ~$C,c>0$~ so that~$$\{z\in\mathbb{C}:c<|z|<C\}\subseteq\sigma_{p}(A)\subseteq\{z\in\mathbb{C}:c\leq|z|\leq C\},$$~ then ~$Ext(A)\subseteq\mathbb{T}$~.
\end{Lemma}
\par This lemma is the key to get the result, we could find the proof in[14.3.3]. With the help of those lemma we could start the most important part:
\begin{Theorem}\label{Theorem15}
{If~$$\varphi(z)=\frac{z+r}{1+rz},0<r<1,$$~~$\varphi$~is a hyperbolic automorphism, then~$Ext(C_\varphi)=\mathbb{T}.$~}\end{Theorem}
{\bf Proof.} First we notice that ~$\mathbb{T}=\{\gamma(it):t\in\mathbb{R}\}$~. Furthermore, ~$\gamma(it)$~ is an eigenvalue of ~$C_\varphi$~ and the correspanding eigenfuction ~$e_{it}$~is bounded on the unit disk. By lemma 2.6 we get ~$\gamma(it)\in Ext(C_\varphi)$~,  then~$\mathbb{T}=\{\gamma(it):t\in\mathbb{R}\}\subseteq Ext(C_\varphi)$~. In other word, ~$Ext(C_\varphi)$~ has rich point spectrum and satisfies the preceding lemma, ~$Ext(C_\varphi)\subseteq\mathbb{T}$~. Hence we complete the proof.
$\hfill\square$

\subsection{hyperbolic nonautomorphism type one}
\par In this section,~$\varphi$~is a hyperbolic, non automorphic linear fractional transformation of the unit disk. First, we consider the class HNA I, which corresponds to the following fixed point configuration: the first fixed point on ~$\partial{D}$~ and the second in ~$\overline{\mathbb{D}}\setminus\ {0}$~ . In this case, we assume that the fixed points are ~$z =1$~ and
~$z=\infty$~. This leads to the standard form~$$\varphi(z) = rz+(1-r)$$~, for some ~$0 < r <1.$~
\par 
Deddens [10, Theorem 3 (iv)] proved that the point spectrum of ~$C_\varphi$~  is the punctured disk~$$
	\sigma_{p}(C_\varphi) = \{{\lambda\in\mathbb{C}:0< |\lambda| < r^{-1/2}}\}.$$~
Moreover, he proved that for every~$ {w}\in\mathbb{C}$~with ~$Re({w}) > -1/2$~, the function ~$e_{w}(z) = (1-z)^{w}$~is an eigenfunction of ~$C_\varphi$~ corresponding to the eigenvalue ~$\lambda= r^{w} $~, that is ~$$C_{\varphi}e_{w}(z)=r^{w}e_{w}(z).$$~As same as we detail in part HA.
\begin{Lemma}\label{lemma21}If ~${\varphi}$~is a hyperbolic nonautomorphism type one of the unit disk,
then ~$C_{\varphi}$~has rich point spectrum.
\end{Lemma}
{\bf Proof.}Back to what we have discussed when we want to make sure HA type composition operator is of rich points. The strategy of the proof this part is exactly same. Namely, we go to the operator ~$$A = C_\varphi,$$~ the open half plane ~$$ w:= {w\in\mathbb{C}: Re (w) >-1/2},$$~ the mapping
~$$h : w \rightarrow H^2(\mathbb{D})$$~ defined by ~$$h(w)(z) = (1-z)^{w},$$~  and the function ~$\gamma:w\rightarrow\mathbb{C}$~ defined by the expression ~$\gamma(w) = r^{w}$~ . Then, conditions (i) and (iii) come true. Regarding condition (ii), it suﬃces to show that the linear span of~$ {e_n : n\in \mathbb{N}_0} $~is dense in Bergman space. For every ~$n \in\mathbb{N}_0$~, the polynomials~$ {(1-z)^k : 0 \leq k \leq n}$~ generate the same linear manifold as the monomials~$ {z^k : 0\leq k \leq n}$~, and by Zhu[]we can also get that ~$ {(1-z)^k : 0 \leq k \leq n}$~ is a total set. Finally, the fact that ~$h$~ is an analytic mapping can be derived in the same way as in the previous proof.
$\hfill\square$
\par It was shown in [14, Theorem 3.1] that the shape of the point spectrum determines the extended spectral picture in the following sense.
\begin{Lemma}\label{lemma22}If ~$A\in B(H)$~ has rich point spectrum and ~$\lambda\in	Ext(A)$~, then~$\lambda int \sigma_{p}(A)\subseteq clos \sigma_{p}(A).$~
\end{Lemma}
\par This result leads to the description of the set of all extended eigenvalues of ~$C_{\varphi}$~.
\begin{Theorem}\label{Theorem17}
{If~$\varphi$~ is HNAI,~$\varphi(z)=rz+1-r$~, for~$ 0<r<1$~,~$Ext(C_\varphi)=\overline{\mathbb{D}}\setminus\ {0}$~}
\end{Theorem}
{\bf Proof.}First, all elements as the type ~$ {(1-z)^{w}: Re{w}>0}$~ are bounded analytic functions, also in Bergman space~$L^{2}_a$~, they can be a total set. By lemma 2.6 we get ~$\overline{\mathbb{D}}\setminus\ {0}\subseteq Ext(C_\varphi).$~On the other hand, last lemma told us that~$if\lambda\in	Ext(C_\varphi)$~, ~$\lambda int\sigma_{p}(C_\varphi)\subseteq clos\sigma_{p}(C_\varphi).$~It means all elements in~$Ext(C_\varphi)$~ are in ~$\overline{\mathbb{D}}\setminus\ {0}$~with lemma 3.1.
In a word, if~$\varphi$~ is HNAI,~$\varphi(z)=rz+1-r$~, for~$ 0<r<1$~, ~$Ext(C_\varphi)=\overline{\mathbb{D}}\setminus\ {0}.$~
$\hfill\square$
\subsection{parabolic automorphism}
\begin{Theorem}\label{Theorem18}
{If~$\varphi$~ is PA,~$$\varphi(z)=\frac{(2-a)z+a}{-az+2+a}$$~where~$ a\in\mathbb{C}$~ and~$ Re(a)=0$~,then~$Ext(C_\varphi)=\partial{D}$~}\end{Theorem}
{\bf Proof.} Similarly, the characteristic function corresponding to parabolic self mapping is also an external function, which naturally satisfies the problem.Some results have be proved that if ~$\varphi$~ is PA, ~$$\varphi(z)=\frac{(2-a)z+a}{-az+2+a}$$~where~$ a\in\mathbb{C}$~ and~$ Re(a)=0$~.The family~$$\{ e_{t}(z)=exp(-t\frac{1+z}{1-z})(t),t\geqslant0 \}$$~ holds for the eigenvectors for~$ C_\varphi$~, and the eigenvalues are ~$e^{-at}(t\geqslant0)$~.It could be considered as ~$\partial{D}$~since ~$at$~ is truly a pure imaginary number and ~$e^{-at}(t\geqslant0)$~ is of model 1 whatever it takes.We can prove that  ~$e_{t}$~ and~$ 1/ e_{t}$~ are bounded analytic functions by maximum modulus principle. Furthermore they are on ~$L^{2}_a$~.Then we know~$$\{ e_{t}(z)=exp(-t\frac{1+z}{1-z})(t),t\geqslant0 \}$$~ can span a dense subset of ~$L^{2}_a$~.Using previous lemma we can obtain the theorem we prepared following.$\hfill\square$

\section{Acknowledgments}
Our thesis focus on the extended eigenvalues of compostion operators, finally we get several presentation about it.

\end {document}